
\documentstyle{amsppt}
\magnification=1200

\hoffset=-0.5pc
\vsize=57.2truepc
\hsize=38truepc
\nologo
\spaceskip=.5em plus.25em minus.20em

 \define\anderone{1}
 \define\atibottw{2}
 \define\atiyjone{3}
 \define\atiysemi{4}
 \define\bauesboo{5}
\define\bottsega{6}
 \define\bouskan{7}
\define\browhueb{8}
\define\curtione{9}
\define\donalkro{10}
\define\duponone{11}
\define\dwykanon{12}
\define\eilmactw{13}
\define\afischon{14}
\define\friehabe{15}
  \define\modus{16}
\define\modustwo{17}
    \define\srni{18}
\define\finite{19}
\define\huebjeff{20}
\define\jamesboo{21}
\define\jeffrtwo{22}
\define\jeffrthr{23}
  \define\kanone{24}
  \define\kantwo{25}
  \define\kanthr{26}
  \define\kanfou{27}
\define\karshone{28}
\define\kauflins{29}
\define\maclaboo{30}
\define\maclbotw{31}
 \define\maybook{32}
\define\milnothr{33}
\define\philston{34}
\define\pressega{35}
\define\puppeone{36}
\define\shulmone{37}
\define\steenone{38}
\define\weinstwe{39}
\define\whitehon{40}
\define\whitehtw{41}
\define\wittesix{42}
\noindent
arch-ive/9505005
\topmatter
\title
Extended moduli spaces
and the Kan construction
\endtitle
\author Johannes Huebschmann{\dag}
\endauthor
\affil
Max Planck Institut f\"ur Mathematik
\\
Gottfried Claren Str. 26
\\
D-53 225 BONN
\\
huebschm\@mpim-bonn.mpg.de
\endaffil
\thanks
{{\dag} The author carried out this work in the framework of the
VBAC research group of EUROPROJ.}
\endthanks
\abstract{Let $Y$ be a CW-complex with a single 0-cell, let $K$ be its Kan
group, a free simplicial group whose realization is a model for the space
$\Omega Y$ of based loops on $Y$, and let $G$ be a Lie group, not necessarily
connected. By means of simplicial techniques involving fundamental results of
{\smc Kan's} and the standard $W$- and bar constructions, we obtain a weak
$G$-equivariant homotopy equivalence from the geometric realization
$|\roman{Hom}(K,G)|$ of the cosimplicial manifold $\roman{Hom}(K,G)$ of
homomorphisms from $K$ to $G$ to the space $\roman{Map}^o(Y,BG)$ of based maps
from $Y$ to the classifying space $BG$ of $G$ where $G$ acts on $BG$ by
conjugation. Thus when $Y$ is a smooth manifold, the universal bundle on $BG$
being endowed with a universal connection, the space $|\roman{Hom}(K,G)|$ may
be viewed as a model for the space of based gauge equivalence classes of
connections on $Y$ for all topological types of $G$-bundles on $Y$ thereby
yielding a rigorous approach to lattice gauge theory; this is illustrated in
low dimensions.}
\endabstract
\date{May 23, 1995} \enddate
\endkeywords
\endsubjclass

\endtopmatter
\document
\leftheadtext{Johannes Huebschmann}
\rightheadtext{Extended moduli spaces and the Kan construction}

\beginsection Introduction

In gauge
theory, one usually studies the space
of gauge equivalence classes of connections
on a principal
bundle or suitable subspaces thereof.
The geometry of the space of connections is quite simple
since it is an affine space.
However its analysis is more intricate,
and
suitable choices of topologies and of completions
must be made, depending on the concrete problem
under consideration.
The miracle is that these analytical problems
disappear
on the space of gauge equivalence classes
of connections.
The present paper
and its successor
\cite\finite \
provide a step towards an
explanation for this.
Usual gauge theory
could be
viewed as non-abelian
{\it singular\/} cohomology
and, in a sense,
we offer here
a corresponding
{\it cellular\/}
approach:
Let $Y$ be a finite CW-complex
with a single 0-cell and $G$ a Lie group,
not necessarily connected.
Let $K$ be the Kan group on $Y$ \cite\kantwo;
this is a simplicial group whose realization
is a model for the  space $\Omega Y$
of based loops on $Y$.
By means of simplicial techniques
involving
fundamental results of {\smc Kan's}
\cite\kantwo\
and the standard $W$- and bar constructions,
we shall obtain
a $G$-equivariant map $\Phi$
from
the realization
$|\Cal H|$
of the
cosimplicial
$G$-manifold
$\Cal H = \roman{Hom}(K,G)$
to
the space
of based maps
from
$Y$ to $BG$ where $G$ acts on $BG$ by conjugation,
and our
main result, Theorem 1.7 below,
will say that
$\Phi$
is a weak homotopy equivalence.
One could say
the domain of $\Phi$
gives a complete set of combinatorial data
which determine a bundle
with a based gauge equivalence class of
connections;
the latter is given by the value
of the data in $\roman{Map}^o(Y,BG)$ under $\Phi$.
We do not know whether
$\Phi$ is in general
a genuine homotopy eqivalence.
For a closed
topological surface,
in Section 2 below, we briefly indicate a construction of
a homotopy inverse of $\Phi$.
When $Y$ is a sphere $S^q,\, q \geq 1$,
with the usual CW-decomposition with only two cells,
the map $\Phi$
boils down to the standard relationship between
$\roman{Map}^o(S^{q-1},G)$
and
$\roman{Map}^o(S^q,BG)\cong \roman{Map}^o(S^{q-1},\Omega BG)$
induced by the
standard map
from $G$ to
the space $\Omega BG$ of based loops on $BG$.
In general,
every topological type of principal $G$-bundle
on $Y$
gives rise to a group of based gauge transformations;
topologically,
the space
$|\Cal H|$
amounts to the union
of the classifying spaces for these groups,
one such space for each topological type.
Here the universal $G$-bundle
$EG \to BG$ is understood
endowed with the universal connection,
as exploited by {\smc Shulman} in his thesis
\cite\shulmone, see
also
our follow up paper
\cite\finite,
and hence a based \lq\lq smooth\rq\rq\  map
from
$Y$ to $BG$
determines a based gauge equivalence
class of connections on its induced bundle.
For a $k$-sphere,
the space of based maps from
a
$(k-1)$-sphere to
$G$
has already been
taken as a model for
the space of based gauge equivalence classes of $G$-connections
on
the $k$-sphere
at various places in the literature,
cf. e.~g. \cite\atiyjone\ (2.3).
Our
construction
offers a generalization thereof,
to
arbitrary (finite) CW-complexes
$Y$, where it yields a
kind of gauge theory on $Y$.
A space similar to
$|\Cal H|$ has been
studied in \cite\afischon.
\smallskip
Why do we resort to
the space $|\Cal H|$
at all?
Apart from the lack of good comparison
between
$\roman{Hom}(\Omega Y,G)$
and
$\roman{Map}^o(Y,BG)$,
for our purposes, there is quite a different reason:
The object
$\Cal H =\roman{Hom}(K,G)$
is
a smooth cosimplicial
manifold
which is {\it finite dimensional\/}
in each degree;
we do not see
how
this could be manufactured
directly from
$\Omega Y$.
This finite dimensionality of
$\Cal H$ in each degree
will be crucial in
\cite\finite:
in that paper we shall carry out
a purely
finite dimensional construction
of the generators of the real cohomology
of $\roman{Map}^o(Y,BG)$
and hence
of the generators of the real cohomology of the offspring moduli spaces
from which
for example
Donaldson polynomials are obtained
by evaluation against suitable fundamental classes
corresponding to moduli spaces of ASD connections.
Another application in \cite\finite\
will be a purely finite dimensional construction
of the Chern-Simons function
for an arbitrary 3-manifold.
This answers a question
raised by {\smc Atiyah}
in \cite\atiysemi\
where he comments
on a possible approach to
the path integral quantization
of the Chern-Simons function.
In fact,
our  paper
\cite\finite\
may be viewed as a step towards
a combinatorial construction of
\lq\lq topological field theories\rq\rq.
Perhaps a suitable quantization thereof then yields
3-manifold invariants of the
Witten-Reshetikhin-Turaev kind,
cf. e.~g. \cite\kauflins.
This would provide a rigorous
construction of Witten's topological quantum field theory
\cite \wittesix.
Our  paper
\cite\finite\
generalizes prior constructions
in
\cite\karshone\
and \cite\weinstwe\
and, furthermore, the subsequent
extensions
thereof
in
\cite\modus,
\cite\modustwo,
\cite\huebjeff,
\cite\jeffrtwo,
\cite\jeffrthr;
in fact, it yields the \lq\lq grand unified theory\rq\rq\
for a {\it general\/} bundle on an {\it arbitrary\/}
compact smooth finite dimensional manifold
searched for by A. Weinstein \cite\weinstwe\
and established by L. Jeffrey \cite\jeffrthr\
over a
surface.
\smallskip
Trying to generalize
the
extended moduli spaces constructed in
\cite\modus,
\cite\modustwo,
\cite\huebjeff,
\cite\jeffrtwo\
to arbitrary bundles
over arbitrary smooth manifolds,
we discovered that
these extended moduli spaces
may be found as suitable subspaces
of
the realization
of
the requisite
cosimplicial  manifolds;
see Section 1 below for details.
This suggests that
the
searched for general extended moduli spaces
should be found within the
realizations of cosimplicial manifolds
of the  aformentioned kind.
We illustrate this in Sections 2 -- 4 below.
\smallskip
In a sense,
the extended moduli space constructions
carried out in
the cited references rely
on the fact that a closed topological surface
different from the 2-sphere
has a combinatorial model which
can entirely be described in terms of the fundamental group
since such a surface
is an Eilenberg-Mac Lane space.
Now for a bundle on an {\it arbitrary\/} space $Y$,
such a naive approach will fail
when
$Y$ is {\it not\/}
an Eilenberg-Mac Lane space.
Our {\it principal innovation\/}
is to
take as combinatorial model
for $Y$
the {\it simplicial nerve\/} (or bar construction) of
the Kan group
$K$ on $Y$.
This idea is behind the constructions
of the present paper, and the structure
of the simplicial nerve
of $K$ will explicitly be exploited
in our follow up paper \cite\finite.
In a sense, the cosimplicial manifold of
homomorphisms from $K$
to the structure group
$G$ generalizes the
usual description of based
gauge equivalence classes of flat connections
in terms of their holonomies
to arbitrary connections.
This statement can be made much more precise:
The geometric realization $|K|$
of $K$
is a topological group
and
the geometric realization
of the cosimplicial manifold
$\Cal H$
of
homomorphisms from $K$
to $G$ amounts to the space
of continuous
homomorphisms from
$|K|$ to $G$.
In the context of smooth bundles
this may look a bit odd
at first
and it seems difficult to view
$|K|$ as a Lie group
but a replacement for
a missing space of
smooth maps from $|K|$ to $G$ is provided
by
what we call the
{\it smooth\/}
geometric realization
$|\Cal H|_{\roman{smooth}}$
of the cosimplicial manifold
$\Cal H$,
cf. Section 1 below;
it is (weakly) homotopy equivalent to
$|\Cal H|$ and may be viewed
as a model for the space
of based gauge equivalence classes of connections.
The lack of decent smooth
structure on $|K|$ is not a problem of infinite dimensions;
artificially,
$|K|$ can be endowed
with a kind of smooth structure
by adjointness
but for our problem of study there is no need to do so
and we do not know what kind of insight
such a smooth structure on $|K|$ would provide.
Our ultimate hope is that framed moduli spaces for various
situations may be found within spaces
of the kind
$|\Cal H|_{\roman{smooth}}$.
\smallskip
For the case of a bundle on a closed surface
$\Sigma$,
the present more general
construction involving a model for the full loop space
rather than a presentation of the fundamental group
of the surface \cite\modus, \cite\modustwo, \cite\huebjeff,
\cite\jeffrtwo, \cite\jeffrthr\
already goes beyond
the earlier extended moduli space
constructions:
The realization
$|\Cal H|$ of
$\Cal H = \roman{Hom}(K\Sigma,G)$
contains the spaces
of based gauge equivalence classes
of {\it all\/}
central Yang-Mills connections
\cite\atibottw,
not just those which correspond
to the absolute minimum or, equivalently,
to projective representations
of the fundamental group $\pi$ of $\Sigma$,
and hence the space
$|\Cal H|$
comes with a kind of Harder-Narasimhan
filtration.
The latter cannot
be obtained from the earlier extended moduli space
constructions. See Section 2 below for details.
Perhaps information about the multiplicative
structure of the cohomology of moduli spaces
can be derived from the resulting models in \cite\finite.
\smallskip
By means of the
simplicial groupoid
constructed in
\cite\dwykanon\
for an arbitrary connected
simplicial set
the present approach
can be extended
to arbitrary
connected simplicial complexes,
in particular,
to triangulated smooth manifolds.
In the non-abelian cohomology spirit,
this will amount
to a {\it simplicial\/} gauge theory.
It may be viewed as a rigorous
mathematics
approach
to lattice gauge theory.
In fact,
the
above cosimplicial manifold
$\Cal H = \roman{Hom}(K,G)$
may be viewed as
a space of {\it parallel
transport functions\/},
cf. e.~g. \cite\philston\
for this notion.
More naively,
given an ordered simplicial complex, viewed as a simplicial set,
contracting a maximal tree yields
a simplicial set with a single
vertex, and the construction
of Kan group can be applied.
To keep the present
paper to size,
we plan to give the details elsewhere.
See also the remark at the end of Section 1.
Our approach
somewhat establishes a link between classical algebraic
topology and the more recent gauge theory
developments in low dimensional topology:
our models for the space of
gauge equivalence classes of connections
involve
classical low dimensional topology notions
such as {\it identity among relations\/}
(Section 3 below)
and
{\it universal quadratic group\/}
(Section 4 below).
\smallskip
Any unexplained notation is the same
as
that in our paper \cite\modus.
Details about cosimplicial spaces
may be found in
 \cite\bottsega\
and \cite\bouskan.
Topological spaces will be assumed endowed
with the compactly generated topology throughout.
\smallskip
I am indebted
to S. Bauer, H. J. Baues, and D. Puppe for discussions,
and to J. Stasheff for a number of most helpful
comments on a draft of the manuscript.
The paper has been written during a stay
at the Max Planck Institut at Bonn.
I wish to express my
gratitude to it
and to its director Professor
F. Hirzebruch for hospitality and support.

\medskip\noindent{\bf 1. The finite model}\smallskip\noindent
Write
$\Delta$ for the category  of finite ordered sets
$[q] = (0,1,\dots,q),\, q\geq 0,$
and monotone maps.
We recall the standard {\it coface\/} and {\it codegeneracy\/} operators
$$
\align
\varepsilon^j
&\colon [q-1] @>>> [q],
\quad
(0,1,\dots, j-1,j,\dots,q-1)
\mapsto (0,1,\dots, j-1,j+1,\dots,q),
\\
\eta^j
&\colon [q+1] @>>> [q],
\quad
(0,1,\dots, j-1,j,\dots,q+1)
\mapsto (0,1,\dots, j,j,\dots,q),
\endalign
$$
respectively.
As usual,
for a simplicial object,
the corresponding face and degeneracy operators
will be written $d_j$ and $s_j$.
Recall that a {\it cosimplicial\/}
object in a category $\Cal C$
is a covariant functor
from $\Delta$
to $\Cal C$.
For example,
the assignment
to
$[q]$ of the standard simplex
$\nabla [q]=\Delta_q$
yields a cosimplicial
space $\nabla$;
here we wish to distinguish
clearly
in notation between
the cosimplicial space
$\nabla$
and the category $\Delta$.
Let $K$ be a free simplicial groupoid,
for example a free simplicial group.
The simplicial structure of $K$
induces
a structure
of {\it cosimplicial\/}
manifold on the groupoid homomorphisms
$
\roman{Hom}(K,G)
$
from $K$ to $G$.
For $q \geq 0$, we shall henceforth write
$
\roman H_q =
\roman{Hom}(K_q,G)
$
so that
$\Cal H=\roman{Hom}(K,G)$
may be depicted as
$
\{\roman H_0,
\roman H_1,
\dots,
\roman H_q,
\dots \}
$
with the requisite smooth maps between
the components induced by monotone maps between finite sets.
\smallskip
The {\it geometric  realization\/}
$
|\Cal H|
$
of $\Cal H$,
cf.
\cite\bottsega,
\cite\bouskan,
is the
space
$
|\Cal H| = \roman{Hom}_{\Delta}(\nabla, \Cal H);
$
this is the subspace
of the infinite product
$$
\roman H_0 \times \roman{Map}(\Delta_1,\roman H_1) \times \dots \times
\roman{Map}(\Delta_q,\roman H_q)
\times \dots
\tag1.1
$$
consisting of all sequences
$
(\phi_0,\phi_1,\dots,\phi_q,\dots)
$
having the property that, for each
monotone map $\theta\colon [i] \to [j]$,
the diagram
$$
\CD
\Delta_i
@>{\theta_*}>>
\Delta_j
\\
@V{\phi_i}VV
@V{\phi_j}VV
\\
\roman H_i
@>{\theta_*}>>
\roman H_j
\endCD
\tag1.2
$$
commutes.
\smallskip
When
$K$ is countable
the geometric realization $|K|$ of $K$
is a
topological groupoid,
cf. e.~g. \cite\milnothr\
where this is proved for simplicial groups.
In general,
one has to take  compactly generated refinements
of the product topologies
on the spaces where
compositions are defined.
Henceforth we suppose $K$ countable.
Then the cosimplicial manifold
$\Cal H=\roman{Hom}(K,G)$
provides a model of
the space
$\roman{Hom}(|K|,G)$
of continuous homomorphisms from
$|K|$ to $G$.
In fact,
for $q\geq 0$,
adjointness
yields a canonical map
from
$
\roman{Map}(\Delta_q,\roman{Hom}(K_q,G))
$
to
$\roman{Map}(K_q \times \Delta_q,G)$,
by construction,
the space
$\roman{Hom}(|K|,G)$
canonically embeds into
the infinite product of the spaces
$\roman{Map}(K_q \times \Delta_q,G)$,
and we have the following tautology:

\proclaim{Proposition 1.3}
Adjointness induces a homeomorphism
between
$\roman{Hom}(|K|,G)$
and
\linebreak
$|\roman{Hom}(K,G)|$.
\endproclaim

More formally,
the geometric realization
$|K|$ is
the {\it coend\/}
$K \otimes _\Delta\nabla$,
cf. e.~g. \cite\maclbotw,
and
we have an adjointness
$$
|\roman{Hom}(K,G)| =
\roman{Hom}_\Delta(\nabla,\roman{Hom}(K,G))
@>>>
\roman{Hom}(K \otimes _\Delta\nabla,G)
=
\roman{Hom}(|K|,G).
$$

\smallskip
Henceforth
we shall exclusively deal with free simplicial groups.
We recall \cite\kantwo\
that a graded set
$X = \{X_0,X_1,\dots\}$, where $X_q \subseteq K_q$, for $q \geq 0$,
is called a set of {\it (free) generators\/} for $K$
provided
$K$ is freely generated by $X$ as a simplicial group.
That is to say:
\roster
\item
If $q \geq 1$ and $0 \leq j <q$ then
$\partial_j x = e_{q-1}$, the neutral element, for every $x \in X_q$.
\item For each $q$,
the set $X_q$ together with all the degeneracies
$s_u s_v \dots s_w x \in K_q$, for $x$ in some
$X_{r}$,
freely generates $K_q$ (as a free group).
\endroster
A set $X$ of free generators together with all its degeneracies
is then called a CW-{\it basis\/} for $K$, and for every $q\leq 1$
and every $x \in X_q$, the value
$\partial_q x \in K_{q-1}$
is called the {\it attaching element\/}
of $x$.
\smallskip\noindent
{\smc Remark.}
Here we give preferred treatment to the {\it last\/} face operator, as is
done in \cite\kanone, \cite\kantwo.
This turns out to be the appropriate thing to do
for principal bundles with  structure group acting
on the {\it right\/} of the total space.
\smallskip
It is proved in \cite{\kantwo\ (2.2)}
that every free simplicial group
has a CW-basis.
By means of a CW-basis, the geometric realization
$|\Cal H|$ of
$\Cal H$ may be cut to size:
In {\smc Anderson's} terminology
\cite\anderone,
the cosimplicial space
$\Cal H$ is {\it primitive\/}
over the projection maps
$p_q$
from
$
\roman H_q= \roman{Hom}(K_q,G)
$
to
$P_q= G^{X_q}$;
this means that, if $\alpha$ runs
over the $\binom qk$ surjections
from $\Delta_q$ to $\Delta_k$
for $k < q$,
the product of $p_q$ and the $p_k \Cal H(\alpha)$
provides a homeomorphism
$$
\roman H_q
@>>>
P_0 \times P_1^{\binom q1} \times \dots\times
P_{q-1}^{\binom q{q-1}} \times P_q.
$$
Given
${
(\phi_0,\phi_1,\dots,\phi_q,\dots)
}$
in
$|\Cal H|$,
for $q \geq 0$,
write $\psi_q\colon \Delta_q \to P_q$
for the composite of $\phi_q$
with the projection
from $\roman H_q$
onto
$P_q$.
For $q \geq 1$,
the \lq\lq last coface map\rq\rq\
$
\varepsilon^q
$
from
$[q-1]$
to
$[q]$
induces
the affine map
from
$\Delta_{q-1}$ to
$\Delta_{q}$
which identifies
$\Delta_{q-1}$ with the last face of $\Delta_{q}$,
that is, with the face opposite the last vertex.
We now consider the product
$$
G^{X_0} \times \roman{Map}(\Delta_1,G^{X_1}) \times \dots \times
\roman{Map}(\Delta_q,G^{X_q})
\times \dots \quad .
\tag1.4
$$
It is finite when
$Y$ is compact.
Henceforth we write
$G^{X_{q}} = e$ when $X_q$ is empty.

\proclaim{Lemma 1.5}
The assignment to
${
(\phi_0,\phi_1,\dots,\phi_q,\dots)
}$
of
${
(\psi_0,\psi_1,\dots,\psi_q,\dots)
}$
induces a homeomorphism
from
$|\Cal H|$
onto
the subspace
$|\Cal H|'$
of {\rm (1.4)}
consisting of all sequences
${
(\psi_0,\psi_1,\dots,\psi_q,\dots)
}$
of maps $\psi_q$
whose restriction
to all but the last faces of $\Delta_q$
has constant value $e \in G^{X_q}$
and which
satisfy the
recursive requirement that,
for each $q$,
the  diagram
$$
\CD
\Delta_{q-1}
@>{(\varepsilon^q)_*}>>
\Delta_{q}
\\
@V{\phi_{q-1}}VV
@V{\psi_{q}}VV
\\
\roman H_q
@>{(\varepsilon^q)_*}>>
G^{X_{q}}
\endCD
$$
commute.
\endproclaim

\demo{Proof}
For $k<q$,
each (affine) surjection from
$\Delta_q$ to $\Delta_k$
induces a continuous map
from
$\roman{Map}(\Delta_k,P_k)$
to
$\roman{Map}(\Delta_q,P_k)$
and these assemble to a continuous map
from
$\roman{Map}(\Delta_k,P_k)$
to
$\roman{Map}(\Delta_q,P_k^{\binom qk})$.
These maps, in turn, assemble
to a continuous map
from $|\Cal H|'$
into (1.1)
which yields a continuous inverse
of the map from
$|\Cal H|$
to
$|\Cal H|'$. \qed\enddemo
\smallskip
Following
\cite\kantwo\
we shall say that a  CW-complex $Y$
is  {\it reduced\/}
provided it has a single 0-cell
and,
for every $(q+1)$-cell $c$,
the characteristic map
$\sigma_c$
from $\Delta_{q+1}$ to $Y$
has values different from
the base point at most on the next to the last face,
that is, on the one opposite to the vertex $A_{q}$ where the
vertices of
$\Delta_{q+1}$
are numbered $A_0,\dots, A_{q+1}$.
We note
that it is uncommon to have
a CW-complex with cells which are images of simplices
but the present description
is an important
ingredient
for
{\smc Kan's}
results which we shall subsequently use.
A
twisting function $t$
from
the first Eilenberg subcomplex
$SY$
of the total
singular complex
of $Y$
to a simplicial group
$K$
is said to be {\it regular\/}
provided (i) the elements
$t(\sigma_c)$
where $c$ runs through the cells of $Y$
of dimension at least one
form the generators of a CW-basis of $K$
and
(ii)
for every subcomplex $Z$ of $Y$,
the image $t(SZ)$
of its first Eilenberg subcomplex
is contained in the simplicial subgroup
of $K$ generated by the
$t(\sigma_c)$
for $c$ in $Z$.
To any reduced CW-complex $Y$, Kan's construction
\cite\kantwo\
assigns
a free simplicial group $KY$ together
with a regular twisting function
$t$ from
$SY$
to $KY$ \cite\kantwo\ and,
furthermore,
to any free simplicial group $K$,
the reverse construction
of
Kan's
assigns a reduced CW-complex
$YK$ together with a regular twisting function
$t$ from $SYK$ to $K$, so that
$YKY \cong Y$ and $KYK \cong K$.
For each
$(q+1)$-cell $c$ with characteristic map $\sigma_c$,
since
$d_{j}\sigma_c$ is the base point
when $j\ne q$,
the twisting function $t$ satisfies
$$
\align
d_i(t \sigma_c) &= t (d_i \sigma_c)= e,\quad 0 \leq i < q,
\\
d_{q}(t \sigma_c) = t (d_{q} \sigma_c) t(d_{q+1}\sigma_c)^{-1}
&= t (d_{q} \sigma_c) \in K_{q-1},
\\
s_i(t \sigma_c) &= t(s_i \sigma_c),\quad
\quad 0 \leq i \leq q,
\\
e_{q+1} &= t(s_{q+1} \sigma_c).
\endalign
$$
See \cite\kantwo\  for details.
The cosimplicial structure of $\Cal H$
may now be described as follows:
For each
$(q+1)$-cell $c$, with characteristic map $\sigma_c$,
write
$G_c$ for the factor of
$\roman H_q = \roman{Hom}(K_q,G)$
which corresponds
to the free generator $t(\sigma_c)$ of $K_q$.
For $0 \leq j < q$,
the composite of the coface map
$\varepsilon^j$
from
$\roman H_{q-1}$ to
$\roman H_{q}$
with the projection onto
$G_c$ is trivial while
the composite
$$
\roman H_{q-1}=
\roman{Hom}(K_{q-1},G)
@>>>
G_c
$$
of the coface map
$\varepsilon^q$
from
$\roman H_{q-1}$ to
$\roman H_{q}$
with the projection onto
$G_c$ is
given by
the assignment to
$\alpha \in
\roman{Hom}(K_{q-1},G)
$
of the value $\alpha (t(d_q(\sigma_c)))$.
The rest of the structure
is now completely determined by the requirement
that $\Cal H$ be a cosimplicial space.
\smallskip
The regularity of the twisting
function $t$
entails  that the total complex of the associated
simplicial
principal bundle
$\pi \colon SY \times _t KY \to SY$
is contractible whence
$KY$ is a {\it loop complex of\/} $SY$
under $t$.
We now explain what this means for us:
The
geometric realization
of $\pi$ is a principal $|K|$-bundle with base $|SY|$.
Pick a homotopy inverse
$\sigma$
from $Y$ to
$|SY|$
of the
counit
$\varepsilon \colon |SY| \to Y$
of the adjointness between the realization and singular complex
functors.
When $Y$ is itself the realization of a simplicial set
there is a canonical such map $\sigma$.
Whether or not this happens to be the case,
$\sigma$ induces a principal $|K|$-bundle
$\kappa \colon P \to Y$ on $Y$
with contractible total space $P$.
In particular,
a standard homotopy theory construction
yields a homomorphism
from the (Moore) loop space
$\Omega Y$
to the geometric realization
$|KY|$
which is a homotopy equivalence.
The bundle $\kappa$
is classified by a map $\rho$ from
$Y$ to
$B|K|$,
and the assignment to
$\phi \in
|\roman{Hom}(K,G)| \cong \roman{Hom}(|K|,G)$
of the composite
$(B\phi) \rho$
yields a $G$-equivariant map
$$
\Phi
\colon
|\roman{Hom}(K,G)|
@>>>
\roman{Map}^o(Y,BG)
\tag1.6
$$
where $G$ acts on $BG$ by conjugation.
By construction,
this map assigns to
$\phi$ a classifying map of the
principal $G$-bundle on $Y$ arising from
the principal $|K|$-bundle $\pi$ via $\phi$.
Notice when
$G$ is discrete,
the space $|\roman{Hom}(K,G)|$
boils down to
the discrete space $\roman{Hom}(\pi_1(Y),G)$
and (1.6)
picks the connected components
of $\roman{Map}^o(Y,BG)$
each of which is contractible.
\smallskip
In general,
bundles over a classifying
space $BH$ of an arbitrary topological group
$H$
are {\it not\/} classified by
representations of $H$ in $G$.
Thus the next result
is somewhat surprising and indicates that
the realization  $|K|$
of the Kan group $K$ has somewhat
special features.

\proclaim{Theorem 1.7}
The map $\Phi$ is a weak $G$-equivariant homotopy equivalence.
\endproclaim

The map $\Phi$ relies on various choices;
under suitable circumstances,
the construction can presumably be made more rigid, in the following way:
When $S$ is a simplicial set
with a single vertex,
for example
arising from an ordered simplicial complex
by contraction of a maximal tree,
the usual Kan group construction
yields a free simplicial group
together
with a twisting function
$t$
from $S$ to $K$, and the
latter
determines a morphism
$\overline t\colon S \to \overline WK$
of simplicial sets,
where
$\overline WK$
refers to the reduced $W$-construction
\cite\maybook.
Its realization
$|\overline t|\colon |S| \to |\overline WK|$,
combined
with a certain universal comparison map from $|\overline WK|$
to
$B|K|$,
cf. e.~g. what is said at the end of
the introduction
to {\smc Steenrod's} paper \cite\steenone,
yields
a map
$
\widetilde t
$
from
$S$ to
$B|K|$,
uniquely determined
by $S$ and the choice of maximal tree.
We do not pursue these issues here.

\smallskip
We now begin with the preparations for the proof of Theorem 1.7.
We shall see the theorem comes down to the
canonical map from $G$ to  $\Omega BG$.
Let $q\geq 1$, and
consider the inclusion
of the $(q-1)$-skeleton $Y^{q-1}$ into the $q$-skeleton $Y^q$.
This is a cofibration with cofibre a one point union
$\vee S^q
$
of as many $q$-spheres as $Y$ has $q$-cells.

\proclaim{Lemma 1.8}
The inclusion of the $(q-1)$-skeleton into the $q$-skeleton
induces a
Hurewicz fibration
$$
|\roman{Hom}(K(\vee S^q),G)|
@>>>
|\roman{Hom}(KY^q,G)|
@>>>
|\roman{Hom}(KY^{q-1},G)|
\tag1.8.1
$$
for the geometric realizations.
\endproclaim

For a one-point union
$\vee T_j$,
the Kan group $K (\vee T_j)$
equals the free product
$* KT_j$ of simplicial groups.
When $Y^{q-1}$ is just the base point,
the assertion
thus amounts to
a homeomorphism between
$
|\roman{Hom}(K(\vee_{X_q} S^q),G)|
\cong
|\roman{Hom}(*_{X_q}K S^q,G)|
$
and
$
\times_{X_q} |\roman{Hom}(KS^q,G)|$,
and there is nothing to prove.

\demo{Proof}
Let $q \geq 2$, and suppose
$Y^{q-1}$ is more than the base point.
Consider the cofibration
$
S^{q-1}
@>>>
B^q
@>>>
S^q
$
of regular CW-complexes,
the spheres
$S^{q-1}$
and
$S^q$
having obvious
such CW-decompositions with two cells.
Inspection shows that
(1.8.1) then boils down to the
standard Hurewicz fibration
$$
\roman{Map}^o(S^{q-1},G)
@>>>
\roman{Map}^o(B^{q-1},G)
@>>>
\roman{Map}^o(S^{q-2},G)
$$
with contractible total space.
In general, the fibration (1.8.1) is induced
via the attaching maps for the $q$-cells
of $Y$ from the product
of such fibrations  involving as many copies as $Y$ has
$q$-cells,
as indicated in the commutative diagram
$$
\CD
|\roman{Hom}(K(\vee S^q),G)|
@>{\cong}>>
\times_{X_{q-1}} \roman{Map}^o(S^{q-1},G)
\\
@VVV
@VVV
\\
|\roman{Hom}(KY^q,G)|
@>>>
\times_{X_{q-1}}\roman{Map}^o(B^{q-1},G)
\\
@VVV
@VVV
\\
|\roman{Hom}(KY^{q-1},G)|
@>>>
\times_{X_{q-1}}\roman{Map}^o(S^{q-2},G)
\endCD
\tag1.8.2
$$
whose bottom map is induced by the attaching maps
of the $q$-cells of $Y$. \qed
\enddemo

\demo{Proof of {\rm (1.7)}}
The map
$\Phi$
is compatible
with the CW-structures and hence induces,
for $n \geq 1$,
a commutative diagram
$$
\CD
\roman{Map}^o(\vee S^{n+1},BG)
@<<<
|\roman {Hom}(K (\vee S^{n+1}),G)|
\\
@VVV
@VVV
\\
\roman{Map}^o
(Y^{n+1},BG)
@<{\Phi^{n+1}}<<
|\roman {Hom}(KY^{n+1},G)|
\\
@VVV
@VVV
\\
\roman{Map}^o
(Y^{n},BG)
@<{\Phi^n}<<
|\roman {Hom}(KY^{n},G)|
\endCD
$$
of fibrations.
Since
for a one-point union
the Kan group
equals the free product
 of the Kan groups for the factors,
in degree one,
the map $\Phi^1$
amounts
to a product of copies
of maps
of the kind
$G \to \Omega BG$,
the number of factors being given by the number of 1-cells of $Y$.
Likewise,
on the fibres, the map comes down to
a product of copies
of maps of the kind
$$
|\roman {Hom}(K  S^{n+1},G)|
@>>>
\roman{Map}^o( S^{n+1},BG),
$$
the number of factors being given by the number of $(n+1)$-cells of $Y$.
However,
$\roman{Map}^o(S^{n+1},BG)$
equals
$\roman{Map}^o(SS^n,BG)$
where \lq $S$\rq\ refers to the based suspension operator,
adjointness identifies
$\roman{Map}^o(SS^n,BG)$
with
$\roman{Map}^o(S^n,\Omega BG)$,
and again we are left with
a standard homotopy equivalence
$$
|\roman {Hom}(K  S^{n+1},G)|
=
\roman{Map}^o(S^n,G)
@>>>
\roman{Map}^o(S^n,\Omega BG).
$$
By induction
we can therefore conclude that
$\Phi$
is a weak  homotopy equivalence.
This
proves the assertion. \qed
\enddemo
\smallskip
As usual, we shall say that
a map from $\Delta_n$ to a smooth manifold $M$
is {\it smooth\/}
when it is defined and smooth on a neighborhood of
$\Delta_n$
in the ambient space.
Henceforth we write
$\roman{Smooth}(\cdot,\cdot)$
for spaces of {\it smooth maps}.
We define the promised {\it smooth
realization\/}
by
$$
|\Cal H|_{\roman{smooth}}
=
|\Cal H|
\cap
G^{X_0} \times \roman{Smooth}(\Delta_1,G^{X_1}) \times \dots \times
\roman{Smooth}(\Delta_q,G^{X_q})
\times \dots .
$$
It is weakly homotopy equivalent to
$|\Cal H|$
and may be viewed as a model for the space
of based gauge equivalence classes of connections
on $Y$ when the latter is a smooth manifold.
\smallskip
Finally
we
explain briefly the notion of attaching element:
We recall \cite\curtione\ that the homotopy groups of $K$
may be described as the homology groups of the {\smc Moore\/}
complex
$$
MK \colon M_0 @<{d_1}<< M_1 @<{d_2}<< M_2 @<{d_3}<< \dots
$$
of $K$.
Here, for $k \geq 1$,
the group $M_k$ is the intersection $\cap_{j=0}^{k-1} \roman{ker}(d_j)$
and the operator $d_k$ in the Moore complex
is the restriction of the last face operator
(denoted by the same symbol).
It is also customary in the literature
to take the intersection of the kernels of the
last face operators and to take the first face operator
as boundary in the Moore complex.
Let
$t(\sigma_c) \in X_q$
be a free generator corresponding to a $(q+1)$-cell
$c$ of $Y$,
attached via the map $\sigma_c$, restricted
to the boundary of $\Delta_{q+1}$;
the latter represents  an element of
$\pi_q(Y^q)$
which,
under the
standard isomorphisms
between
$\pi_q(Y^q)$,
$\pi_{q-1}(\Omega Y^q)$,
and $\pi_{q-1}(KY^q)$,
passes
to
the class
in
$\pi_{q-1}(KY^q)$
represented
by the value
$\partial_q x \in K_{q-1}$
of the {\it attaching element\/}
$t(\sigma_c)$.
\smallskip\noindent
{\smc Remark.}
In \cite\dwykanon,
the relationship
between
reduced simplicial sets and simplicial groups
has been extended to one between
connected simplicial sets and simplicial groupoids.
By means of it,
we intend to generalize elsewhere
the above constructions to
arbitrary simplicial complexes
and in particular to
triangulated
smooth manifolds.
This will
enable us to remove the seemingly
fuzzy notion of reduced CW-complex which is somewhat unnatural
for smooth manifolds.
More naively,
cf. \cite\kanthr.
given an ordered simplicial complex
$Y$, viewed
as a simplicial set,
contracting a maximal tree $T$ yields
a
simplicial set $Y/T$,
and the Kan construction applied to it
yields a free simplicial group
$K$ and a twisting function
$\widetilde t$
from
$Y/T$ to $K$ so that
$K$ is a loop complex for $Y/T$;
composing with the projection
from $Y$ to
$Y/T$
we obtain a twisting
function
$t$ from
$Y$ to $K$ so that
$K$ is a loop complex for $Y$, that is,
the resulting
simplicial principal bundle
$Y\times _t K \to Y$
has contractible total space.
A similar theory can be made for
$K$
with the modification that
the simplices of $Y$ will not constitute
a CW-basis of $K$.
It remains to be seen
which approach is the most suitable one
for what kind of problem.

\medskip\noindent{\bf 2. Closed surfaces}\smallskip\noindent
Let $\Sigma$ be a closed topological surface
of genus $\ell \geq 0$,
endowed with the usual CW-decomposition
with a single 0-cell $o$,
with 1-cells
$u_1,v_1,\dots,u_\ell,v_\ell$, and with a single 2-cell
$c$.
We suppose the decomposition regular in the above sense.
For $1 \leq j \leq \ell$, write
$x_j$ and $y_j$ for the based homotopy class of
$u_j$ and $v_j$ respectively, and denote by $r$
the based homotopy class of the attaching map for $c$.
Then
$$
\Cal P = \langle
x_1,y_1,\dots,x_\ell,y_\ell; r\rangle
$$
is a presentation for the fundamental group $\pi$ of $\Sigma$.
We suppose things have been arranged in such a way that
$r= \Pi [x_j,y_j]$
in the free group $F$ on the generators.
When the genus is zero, $\Cal P$ is to be interpreted
as a non-trivial presentation of the trivial group,
with $F$ the trivial group.
The Kan group
$K = K\Sigma$ for $\Sigma$ is the free simplicial group
with $K_0 = F$,
with $K_1$ the free group on
$2\ell +1$ generators
$r, s_0(x_1), s_0(y_1),\dots, s_0(x_\ell), s_0(y_\ell)$
where only $r$ is non-degenerate
and,
for $q \geq 2$,
$K_q$
is the free group on the
$(2\ell +q)$
degenerate generators
$$
s_q s_{q-1} \dots s_0(x_j),\quad
s_q s_{q-1} \dots s_0(y_j),\quad
s_{j_q} s_{j_{q-1}} \dots s_{j_1} r,\quad
q \geq j_q > j_{q-1} > \dots > j_1 \geq 0.
$$
Moreover,
the {\it only\/} face operators which are {\it not\/} determined by
the simplicial identities  are
$$
d_0(r) = e,
\quad
d_1(r) = \Pi [x_j,y_j],
$$
and the degeneracy operators are completely determined by the construction
itself.
In particular,
for genus $\ell \geq 1$,
the {\it Moore\/} complex
of $K$
has zero'th homology group
$\pi_1(\Sigma)$
and is exact in higher dimensions.
\smallskip
The relator $r$ induces a smooth  map
from $G^{2\ell}$
to $G$ in the usual way,
where
$G^{2\ell}$
is interpreted to be a single point when $\ell$ is zero;
abusing notation, we denote this map by $r$ as well.
The geometric realization
$|\Cal H|$ of the resulting cosimplicial
space
$\Cal H = \roman{Hom}(K,G)$
is the {\it fibre\/} of $r$.
In fact, the diagram (1.8.2), with $q=2$, now boils down to
$$
\CD
|\roman{Hom}(K(S^2),G)|
@>{\cong}>>
\Omega G
\\
@VVV
@VVV
\\
|\roman{Hom}(K\Sigma,G)|
@>>>
\roman{Map}^o(B^1,G)
\\
@VVV
@VVV
\\
G^{2\ell}
@>{r}>>
G
\endCD
\tag2.1
$$
where $\Omega G$ refers to the space
$\roman{Map}^o(S^1,G)$
of based loops as usual
and $B^1$ to the closed interval.
In particular, when the genus $\ell$ is zero,
$|\roman{Hom}(K\Sigma,G)|$
amounts to $\Omega G$.
This illustrates once more the well known relationship
between moduli spaces over a complex curve and the loop group,
cf. \cite\pressega.
\smallskip
The topological type of the corresponding bundles,
that is, the connected components of the
realization
$|\Cal H|$, may be described as follows:
We take the description of
$|\Cal H|$
as the fibre of $r$, that is to say,
$|\Cal H|$
is now the space of pairs
$(w,\phi)$ where  $w \in G^{2\ell}$ and
$\phi \colon I \to G$ is a path in $G$ from
$e$ to $r(w)$.
Given a point
$(w,\phi)=(w_1,w_2,\dots w_{2\ell -1},w_{2\ell},\phi)$
of
$|\Cal H|$,
pick paths $u_j$ in $G$ from $e$ to $w_j$
and let $\psi \colon I \to G$
be the path in $G$ from $e$ to $r(w)$
given by
$\psi(t) = [u_1(t),u_2(t)] \dots [u_{2\ell -1}(t),u_{2\ell}(t)]$.
Then the composite $\psi^{-1} + \phi$
is a closed path in $G$ from $e$ to $e$;
its class in $\pi_1(G)$
represent the topological type or
connected component of
$|\Cal H|$ in which $(w,\phi)$ lies.
\smallskip
We now show how the
based gauge equivalence classes of the
critical sets of the Yang-Mills
functional
for the gauge theory
over $\Sigma$ with reference to the group $G$
\cite\atibottw\
can be found in the space
$|\Cal H|$:
Suppose at first that $\Sigma$ is not a 2-sphere.
Write $\pi = \pi_1(\Sigma)$
and view $\pi$ as a simplicial group
$\{\pi_q\}$
with $\pi_q = \pi$ for each $q$ and all face and degeneracy
operators the identity map.
The
canonical projection
of simplicial groups
from $K = K\Sigma$ to $\pi$
has kernel
the Kan group $KY$
where $Y$ arises from the
universal covering $\widetilde \Sigma$
of $\Sigma$
by contraction of a maximal tree to a point,
and there results an extension
$$
1
@>>>
KY
@>>>
K
@>>>
\pi
@>>>
1
$$
of simplicial groups.
Their realizations yield the extension
$$
1
@>>>
|KY|
@>>>
|K|
@>>>
\pi
@>>>
1
$$
of topological groups
and the group $\pi_0|KY|$ of components
may be identified with the kernel $N$ of the projection from
$F = K_0$ to $\pi$; notice it coincides with the fundamental group
of the 1-skeleton of $\widetilde \Sigma$.
Dividing out $[N,F]$
we obtain the groups
 $N \big / [N,F] \cong \bold Z$
and
$\Gamma =  N \big / [N,F]$
which
yield the universal central extension
$$
1
@>>>
\bold Z
@>>>
\Gamma
@>>>
\pi
@>>>
1
$$
of $\pi$,
the central copy $\bold Z$ being generated by $r[N,F]$.
The injection of $\bold Z$ into the reals $\bold R$
then induces the central extension
$$
1
@>>>
\bold R
@>>>
\Gamma_{\bold R}
@>>>
\pi
@>>>
1
$$
of $\pi$.
The projection of
$|KY|$
onto
its group
$\pi_0|KY|$ of components,
combined with
the projection onto $\bold Z  \cong N /[F,N]$,
extends to a continuous homomorphism
$\vartheta$
from
$|KY|$
to
$\bold R$.
In fact, adjointness yields
a bijection
$$
\roman{Hom}(KY,S\bold R)
@>>>
\roman{Hom}(|KY|,\bold R)
$$
where $S\bold R$
is the singular complex of $\bold R$, viewed as a simplicial
group,
and it is straightforward to
extend the assignment of $1$ to
$r$
to a morphism of simplicial groups
from $KY$ to $S \bold R$.
The homomorphism $\vartheta$, in turn, induces a continuous
surjective homomorphism $\Theta$ from $|K|$ to
$\Gamma_{\bold R}$
as indicated in the commutative diagram
$$
\CD
1
@>>>
|KY|
@>>>
|K|
@>>>
\pi
@>>>
1
\\
@.
@V{\vartheta}VV
@V{\Theta}VV
@V{\roman{Id}}VV
@.
\\
1
@>>>
\bold R
@>>>
\Gamma_{\bold R}
@>>>
\pi
@>>>
1
\endCD
$$
of extensions of topological groups.
The homomorphism $\Theta$ induces an injection
$$
\roman{Hom}(\Gamma_{\bold R}, G)
@>>>
\roman{Hom}(|K|, G) = |\Cal H|.
$$
The space
$\roman{Hom}(\Gamma_{\bold R}, G)$
is well known to be
that of based gauge equivalence classes of the
critical sets of the Yang-Mills
functional (for all topological types of bundles)
\cite\atibottw.
Formally,
the subspace
$\roman{Hom}(\Gamma_{\bold R}, G)
$
of
$|\Cal H|$
decomposes the latter
into $G$-equivariant \lq\lq Morse strata\rq\rq;
in fact,
it yields a kind of
Harder-Narasimhan
filtration
of $|\Cal H|$,
and the resulting decomposition
of the latter
is a kind of generalized Birkhoff decomposition,
cf. \cite\pressega\
and what is said below.
There is even an obvious candidate for a Morse
function arising from the energy of the paths in
$\roman{Smooth}^o(B^1,G) \subseteq \roman{Map}^o(B^1,G)$,
cf. (2.1);
note that
$B^1$
is just the unit interval,
and we run into a certain variational problem
with additional boundary constraints
coming from the word map $r$.
Details have not been worked out yet.
\smallskip
We now explain briefly
how
under the present circumstances
a homotopy inverse of the map
$\Phi$,
cf. (1.6) above,
may be obtained.
Given a smooth principal bundle $\xi$
on $\Sigma$,
the holonomy yields a
smooth map
from
the space
$
\Cal A(\xi)
$
of connections to
$\roman{Hom}(|K|, G) = |\Cal H|$
which, after a suitable choice
of $\vartheta$ has been made,
restricts to a map from the space
of Yang-Mills connections
to
$\roman{Hom}(\Gamma_{\bold R}, G)$.
More precisely,
with the present conventions,
the surface $\Sigma$ is obtained from a 2-simplex
$\Delta_2$ with vertices $A_0, A_1, A_2$
in such a way that
its characteristic map
$\sigma$ sends the faces $(A_0,A_1)$ and $(A_1,A_2)$
to the 0-cell
and
the face $(A_0,A_2)$ to the boundary path
$\Pi [u_j,v_j]$.
For each point $p$ of the first face $(A_1,A_2)$,
let $w_p$ be the linear path in $\Delta_2$ joining
the vertex $A_0$ with $p$.
The assignment to a
connection $A$ on $\xi$ of the holonomies
of the closed paths $\sigma \circ u_j$ and
$\sigma \circ v_j$
yields a smooth map
from the space
$\Cal A(\xi)$ of connections
to $G^{2\ell}$,
and the assignment to $A$
of the holonomies of the closed paths
$\sigma \circ w_p$
yields a lift of this map
to the space $|\Cal H |_{\roman{smooth}}$
which is smooth in a suitable sense.
Assembling these maps over all topological types
of bundles
we obtain in fact a homotopy inverse of
the above map $\Phi$.
The existence of this map is due to the fact that we are working
over a topological surface
where the combinatorics of the situation
is simple.
Finer combinatorial
tools will presumably
yield a homotopy inverse of $\Phi$
in general.
\smallskip
When $\Sigma$ is the 2-sphere,
in view of the identification
of
$\roman{Hom}(|K|,G)$ with $\roman{Map}^o(S^1,G)=\Omega G$,
with $G=S^1$, we see there is a surjective homomorphism
from $|K|$ to $S^1$ which classifies the universal cover
of
$|K|$.
For general $G$, this surjection induces an embedding
of $\roman{Hom}(S^1,G)$ into
$\roman{Hom}(|K|,G)=\Omega G$.
The space
$\roman{Hom}(S^1,G)$
is well known to be
that of based gauge equivalence classes of the
critical sets of the Yang-Mills
functional over the 2-sphere \cite\atibottw,
cf. also \cite\friehabe,
yielding the Birkhoff decomposition,
cf. \cite\pressega.
\smallskip
We conclude this Section with a topological remark.
The fibrations
(1.8.1)
are known to be rationally trivial,
at least for $G$ simply connected.
Some hints for the general case may be found in
\cite\donalkro,
and, for the present special case,
where (1.8.1) boils down to the
left-hand vertical fibration in (2.1),
the rational triviality
may be found in \cite\atibottw.
It also follows from
Theorem 7.1 in \cite\finite.
However, over the integers,
the fibration under discussion
is in general certainly not trivial.
We briefly explain this for $Y$ a closed surface $\Sigma$:
Its is well known that the
word map
$r$
from $G^{2\ell}$ to $G$
map is not null homotopic unless $G$ is abelian,
cf. \cite\jamesboo.
For example, for $G=\roman{SU}(2)$,
the commutator map from
$G \times G$ to $G$ factors through
$S^6 = S^3 \wedge S^3$
and generates $\pi_6(S^3)$ which is finite cyclic of order 12.
In fact this generator is the {\smc Samelson\/}
product $[a,a]$ where $a$ refers to the generator
of $\pi_3(S^3)$.
A general word map $r$ produces
higher degree generators
in the homotopy
of $G$ whence $r$ will certainly not
be null homotopic.
In particular, with
coefficients in a finite field,
the spectral sequence
of the fibration will in general be non-trivial.

\medskip\noindent
{\bf 3. 3-complexes and 3-manifolds}
\smallskip\noindent
Let $Y$ be a
3-complex with a single 3-cell,
for example,
a closed compact 3-manifold,
endowed with a regular CW-decomposition
with a single 0-cell $o$,
with 1-cells
$u_1,\dots,u_\ell$, 2-cells
$c_1,\dots,c_\ell$, and a single
3-cell $c$.
For $1 \leq j \leq \ell$, write
$x_j$ and $r_j$ for the based homotopy class of
$u_j$ and $c_j$ respectively, and denote by $\sigma$
the based homotopy class of the attaching map for $c$.
Then
$$
\Cal S = \langle
x_1,\dots,x_\ell; r_1,\dots,r_\ell; \sigma\rangle
$$
is a
{\it spine\/}
for $Y$;
in particular,
(i) t
$
\Cal P
=\langle x_1,\dots,x_\ell; r_1,\dots,r_\ell\rangle
$
constitutes a
presentation of the fundamental group
$\pi$ of
$Y$
so that
the attaching maps
of the 2-cells
assign a word $w_j$
in the free group $F$ on the generators
to each relator $r_j$,
and (ii)  the attaching map
$\sigma$ of the single 3-cell assigns
an
{\it identity among relations\/}
$$
i = z_1 r^{\varepsilon_1}_{j_1}z_1^{-1} \dots
z_m r^{\varepsilon_m}_{j_m}z_m^{-1}
\tag3.1
$$
to $c$
representing
the element of the second homotopy group
$\pi_2(Y^2)$
of the 2-skeleton of
$Y$ which is killed by the 3-cell $c$;
here each $z_k$ is an element of $F$,
and the meaning of \lq\lq identity among relations\rq\rq\
will be made clear below in terms of the structure of the Kan group.
See \cite\browhueb\
for more details on the notion of identity among relations.
\smallskip
To spell out the Kan group
$K = KY$ for $Y$,
we do not distinguish in notation
between
the values
of the characteristic maps
of the cells
under the twisting function
$t$ from $SY$ to $K$
and the based homotopy classes
in the spine $\Cal S$
they
correspond to.
With these preparations out of the way,
the group $K$
is the free simplicial group
with $K_0 = F$,
with $K_1$ the free group on
$2\ell$ generators
$r_1,\dots,
r_\ell, s_0(x_1),\dots, s_0(x_\ell)$,
the $r_1,\dots,
r_\ell$ being
non-degenerate,
with $K_2$ the free group on
$3\ell$ degenerate generators
$$
s_0(r_1),\dots,
s_0(r_\ell),
s_1(r_1),\dots,
s_1(r_\ell),
s_1s_0(x_1),\dots, s_1s_0(x_\ell)
$$
together with a single
non-degenerate generator
$\sigma$,
and,
for $q \geq 3$,
$K_q$
is free on
a certain number of
degenerate generators.
Moreover,
the {\it only\/} face operators which are {\it not\/} determined by
the simplicial identities  are
$$
\align
d_0(r_j) &= e,
\quad
d_1(r_j) = w_j \in K_0,
\quad
d_0(\sigma) = e,
\quad
d_1(\sigma) = e,
\\
d_2(\sigma) &=
(s_0z_1) r^{\varepsilon_1}_{j_1}(s_0z_1)^{-1} \dots
(s_0z_m) r^{\varepsilon_m}_{j_m}(s_0z_m)^{-1} \in K_1,
\quad \varepsilon_j = \pm 1,
\endalign
$$
and the degeneracy operators are completely determined by the
simplicial identities as well.
We note that $i$ to be an identity among relations
means precisely that
$d_1 d_2(\sigma) = e$ or, equivalently, that the
word $i$ in the generators of $F$ arising from substituting
each $w_j$ for $r_j$
reduces to the trivial element of $F$.
In particular,
the
part
$
M_0 @<{d_1}<<
M_1 @<{d_2}<<
M_2
$
of the {\it Moore\/} complex
of $K$
determines an exact sequence
$$
1
@<<<
\pi_1(Y)
@<<<
M_0
@<<<
M_1 \big / (d_2 M_2)
@<<<
\pi_2(Y)
@<<<
1
$$
which is just the part
$$
1
@<<<
\pi_1(Y)
@<<<
\pi_1(Y^1)
@<<<
\pi_2(Y,Y^1)
@<<<
\pi_2(Y)
@<<<
1
$$
of the long exact homotopy sequence
of the pair
$(Y,Y^1)$.
The resulting cosimplicial manifold $\Cal H = \roman{Hom}(K,G)$
has
$
\roman H_0 = G^{\ell},
\
\roman H_1 = G^{2\ell},
\
\roman H_2 = G^{3\ell+1},
$
etc., the coface and codegeneracy maps
being determined
by the simplicial structure of $K$
spelled out above.
\smallskip
The $\ell$-tuple $(r_1,\dots, r_\ell)$ of relators induces a smooth  map
from $G^{\ell}$
to
$G^{\ell}$
in the usual way which
we denote by $r$ with an abuse of notation
where
$G^{\ell}$
is interpreted to be a single point when $\ell$ is zero, and
the geometric realization
of the cosimplicial
space
$\roman{Hom}(KY^2,G)$
is the {\it homotopy fibre\/} of $r$,
as inspection of
the diagram (1.8.2) with $q=2$ shows.
Consequently, the geometric realization
of $\Cal H = \roman{Hom}(K,G)$
admits the following description:
The attaching map of the single 3-cell
of $Y$
induces
a homomorphism
of free simplicial groups
from
$KS^2$
to
$KY^2$
which is given by the assignment
to a free generator
of
the free cyclic group
$K_1(S^2)$
of
$$
(s_0z_1) r^{\varepsilon_1}_{j_1}(s_0z_1)^{-1} \dots
(s_0z_m) r^{\varepsilon_m}_{j_m}(s_0z_m)^{-1} \in K_1 = K_1(Y^2)
= K_1(Y).
$$
This homomorphism induces a map
$\sigma^*$ from
$
|\roman{Hom}(KY^2,G)|
$
to
$
|\roman{Hom}(KS^2,G)|
=
\Omega G,
$
and the realization
$|\roman{Hom}(KY,G)|$
is the homotopy fibre of the map $\sigma^*$.

\medskip\noindent{\bf 4. Simply connected polyhedra and 4-manifolds}
\smallskip\noindent
A simply connected 4-manifold $Y$ may be written as
the cofibre of a map $f$
from
the 3-sphere $S^3$ to a bunch $\vee_{\ell} S_j^2$
of $\ell$ copies of the 2-sphere.
In general this construction yields a simply connected
4-complex with a single 4-cell $c$,
with characteristic map
$\sigma$ from $\Delta_4$ to $Y$; it is of the homotopy type
of a 4-manifold
if and only if the attaching map $f$
induces
a non-degenerate
quadratic form
on
$\roman H_2Y \cong \bold Z^{\ell}$.
However,
for non-degenerate intersection form,
it may {\it not\/} yield all
{\it smooth\/} simply connected 4-manifolds
and finer decompositions might be necessary
to cover these.
Likewise we could model a non-simply connected
4-manifold.
We concentrate here on the present
situation
and
for the moment we work with a general 4-complex
$Y$ arising from an arbitrary attaching map $f$.
The corresponding Kan group $K = KY$
has $K_0$ trivial,
$K_1$ the free group on $\ell$ generators
$t(\sigma_1),\dots,
t(\sigma_{\ell})$,
where
$\sigma_1,\dots,\sigma_{\ell}$
are the characteristic maps of the
2-cells of $Y$,
$K_2$ the free group on the $2\ell$ degenerate generators
$$
s_0t(\sigma_1),\dots,
s_0t(\sigma_{\ell}),
s_1t(\sigma_1),\dots,
s_1t(\sigma_{\ell})
\tag4.1
$$
and
$K_3$ the free group on the $3\ell$ degenerate generators
$$
s_1s_0t(\sigma_1),\dots,s_1s_0t(\sigma_{\ell}),
s_2s_0t(\sigma_1),\dots,s_2s_0t(\sigma_{\ell}),
s_2s_1t(\sigma_1),\dots,s_2s_1t(\sigma_{\ell})
$$
together with a single non-degenerate generator
$t(\sigma)$.
The only face operators which are not determined by the simplicial
identities are
$$
d_0(t\sigma) = d_1(t\sigma) = d_2(t\sigma) = e,
\quad
d_3(t(\sigma) = t(d_3(\sigma)).
$$
Notice $d_3(\sigma)$
is a singular 3-simplex
of $Y$ and
$t(d_3(\sigma)) \in K_2$
is a word in the
free generators (4.1);
in analogy with what was said in previous Sections,
we write
$r =t(d_3(\sigma)) \in K_2$.
\smallskip
We now explain how this element $r$ may be made explicit.
To this end we
recall
that, by a result of J. H. C. Whitehead,
$\pi_3(Y^2) = \pi_3(\vee_{\ell} S_j^2)$
equals the universal quadratic group
$\Gamma(\pi_2(Y))$
\cite\eilmactw\
on $\pi_2(Y)$ and hence
is free abelian of rank $\binom{\ell+1}2$,
cf. \cite{\bauesboo, \whitehon, \whitehtw};
more explicitly,
after a choice
$a_j \in \pi_2(S_j^2)\cong \bold Z$
and
$b_j \in \pi_3(S_j^2)\cong \bold Z$
of  generators
has been made,
$1 \leq j \leq n$,
a basis
of $\pi_3(Y^2)$ is given by
$b_1,\dots, b_\ell$
and the
Whitehead products
$[a_i,a_j]$ for $i<j$.
We now translate this
to the second homotopy group
$\pi_2(KY^2)$  ($\cong \pi_3(Y^2)$ )
of the Kan group on the 2-skeleton $Y^2$:
For a single 2-sphere $S^2$,
the
Kan
group
$KS^2$ has
$K_1$
free cyclic with generator
$x=t(\sigma_1)$
and the Moore complex
$MS^2$
has $M_0 = e$,
$M_1 = \bold Z$,
and
$M_2$ the commutator subgroup of
$K_2 = s_0K_1 * s_1K_1$.
The first homology group
of $MS^2$
equals $K_1$; this
is a copy of the integers as it should be
since it is just $\pi_2 S^2$ and,
likewise,
the second homology group of $MS^2$
is a copy of the integers,
generated by the commutator
$[s_0(x),s_1(x)] \in [K_2,K_2]$;
this element corresponds to the Hopf map from
$S^3$ to $S^2$ which generates $\pi_3(S^2)$.
See \cite \kanthr\ (p. 310) for details.
\smallskip
We now return to our
2-complex
$Y^2= \vee_{\ell} S_j^2$.
For simplicity,
for $1 \leq j \leq \ell$,
write
$x_j = t(\sigma_j) \in K_1$
for the free generators
corresponding to the 2-spheres in $Y$.
The
Kan
group
$KY^2$ has
$
K_2= s_0 K_1 * s_1K_1$
and
$K_3= s_1s_0K_1 * s_2s_0K_1* s_2s_1K_1$
etc.,
and the Moore complex
of
$Y^2$
has
first homology group
the group $K_1$ made abelian and second homology group
generated by the classes of the commutators
$$
v_j
=[s_0(x_j),s_1(x_j)] \in [K_2,K_2],
\quad
1 \leq j \leq \ell,
$$
and of the elements
$$
w_{i,j} = s_0(x_i) v_j (s_0(x_i))^{-1}
\in [K_2,K_2],
\quad
1 \leq i<j \leq \ell.
$$
The elements $v_j$
and $w_{i,j}$
correspond to the
generators written $b_j$ and
$[a_i,a_j]$ above, respectively.
The former assertion is obvious and the latter one
may be seen by inspection of the
long exact homotopy
sequence of the extension
$$
1
@>>>
[K,K]^{\roman{Ab}}
@>>>
K\big/[[K,K],[K,K]]
@>>>
K^{\roman{Ab}}
@>>>
1
$$
of simplicial groups:
In fact,
the canonical maps from
$[K,K]$ to $K$
and
$[K,K]$
to
$[K,K]^{\roman{Ab}}$
induce an isomorphism from
$\pi_2(K)$
onto
$\pi_2[K,K]^{\roman{Ab}}$,
and the action of
$\pi_1K^{\roman{Ab}}$
on
$\pi_2[K,K]^{\roman{Ab}}$
corresponds to the operation of
Whitehead product in
$\pi_2(Y^2)$.
The attaching element
$r$ is now a word in
the $v_j$ and the $w_{i,j}$.
Moreover the quadratic form
on the second integral cohomology
may be described in the following way:
The relevant part of {\smc Whitehead's\/}
exact sequence \cite\whitehtw\
looks like
$$
\roman H_4(Y)
@>b>>
\Gamma(\pi_2(Y))
@>>>
\pi_3(Y)
$$
and the quadratic map
from
$\pi_2(Y) = \roman H_2(Y)$
to
$\roman H_2(Y) \otimes \roman H_2(Y)$
given by
the assignment to $a$ of
$a \otimes a$
factors through a homomorphism of abelian groups
from
$\Gamma(\pi_2(Y))$
to
$\roman H_2(Y) \otimes \roman H_2(Y)$.
The composite with the boundary $b$ yields a homomorphism
from
$\roman H_4(Y)$
to
$\roman H_2(Y) \otimes \roman H_2(Y)$
the dual of which is the intersection pairing
on $Y$. In particular $Y$ is the homotopy type
of a simply connected 4-manifold if and only if
the intersection pairing is non-degenerate.
The non-degeneracy of the intersection pairing now translates
in an obvious way
to a condition on the word map $r$.
Notice the similarity of the situation
with that over a surface,
cf. Section 2 above.
\smallskip
The cosimplicial manifold
$\Cal H = \roman{Hom}(K,G)$
has
$\roman H_0 = e$,
$\roman H_1 = G^{\ell}$,
$\roman H_2 = G^{2\ell}$,
$\roman H_3 = G^{3\ell +1}$,
and the only part of the cosimplicial structure
which is not determined by the structure itself is
the composite
of
$\varepsilon^3
$
from
$\roman H_2$
to
$\roman H_3$
with the projection onto the primitive part
$G$ of
$\roman H_3$
which corresponds to the single 4-cell of $Y$.
We write this as a word map
$
r
$
from
$G^{\ell} \times G^{\ell}$
to $G$.
This makes perfect sense
since it is given by the assignment to
$(c_1,\dots,c_\ell, d_1,\dots,d_\ell) \in G^{2\ell}$
of
the element of $G$ which is obtained
by substituting
$c_j$ and $d_j$
for
each occurrence
of $s_0(x_j)$
and
of $s_1(x_j)$,
respectively,
in
$r =t(d_3(\sigma)) \in K_2$.
The smooth geometric realization of $\Cal H$
may now be described as the space of pairs
$(\phi_1, \phi_3)$ of smooth maps
$\phi_1 \colon \Delta_1 \to
\roman H_1 = G^{\ell}$
and
$\phi_3 \colon \Delta_3 \to G$
subject to the conditions
\roster
\item $\phi_1(0) = \phi_1(1) = e$,
\item $\phi_3$ has constant value $e$ on the first three
faces of $\Delta_3$, and
\item the diagram
$$
\CD
\Delta_2
@>{\varepsilon^3}>>
\Delta_3
\\
@V{(\phi_1 \circ \eta^0,\phi_1 \circ \eta^1)}VV
@VV{\phi_3}V
\\
G^{\ell} \times G^{\ell}
@>>r>
G
\endCD
\tag4.2
$$
\endroster
is commutative.
In some more detail,
realize the standard simplex
$\Delta_q$ in $\bold R^{q+1}$
as usual as the subset
of points $(t_0,\dots,t_q)$
defined by
$t_j \geq 0$ and $\sum t_j =1$
so that the maps
$\eta^0$ and
$\eta^1$
from
$\Delta_2$ to
$\Delta_1$
are given by
$$
\eta^0(t_0,t_1,t_2)
=
(t_0+t_1,t_2),
\quad
\eta^1(t_0,t_1,t_2)
=
(t_0,t_1+t_2).
$$
Notice
the resulting map $(\eta^0,\eta^1)$
from
$\Delta_2$
to
$\Delta_1 \times \Delta_1$
identifies
$\Delta_2$
with one of the two simplices
in the triangulation
of
$\Delta_1 \times \Delta_1$
coming into play
in the {\it shuffle map\/},
cf. p. 243 of \cite\maclaboo.
When we take
$(t_1,\dots,t_q)$
as independent variables
on $\Delta_q$,
the realization
of $\Cal H$
appears as the space of
pairs of $G$-valued smooth maps $(\phi_1,\phi_3)$,
where
$\phi_1$ is a smooth function of a single variable
$t\in I$
while
$\phi_3$ is a smooth function of three variables
$t_1,t_2,t_3$ defined for
$t_1+t_2+t_3 \leq 1$
and
$t_j \geq 0$,
subject to the conditions
\roster
\item
$\phi_1(0) = \phi_1(1) = e$,
\item
$\phi_3(t_1,t_2,t_3) = e$ if $t_1 = 0$, $t_2 = 0$, or if
$t_1+t_2+t_3 = 1$, and
\item
$\phi_3(t_1,t_2,0) = r(\phi_1(t_1),\phi_1(t_1+t_2))$.
\endroster
When $Y$ underlies a smooth
4-manifold, this space of maps is a model for the space
of based gauge equivalence classes of
connections on all topological types of bundles
on $Y$.
Perhaps moduli spaces
of
based gauge equivalence classes of
ASD-connections
can be found
within this space.
\smallskip
We conclude this Section with a remark on the topology
of the space of based gauge equivalence classes
of all connections:
Under the present circumstances,
$|\roman{Hom}(KY^2,G)|
=
\times_{\ell}\Omega G$,
and the diagram (1.8.2) boils down to
$$
\CD
|\roman{Hom}(K(S^4),G)|
@>{\cong}>>
\roman{Map}^o(S^3,G)
\\
@VVV
@VVV
\\
|\roman{Hom}(KY,G)|
@>>>
\roman{Map}^o(B^3,G)
\\
@VVV
@VVV
\\
\times_{\ell}\Omega G
@>{\tau}>>
\roman{Map}^o(S^2,G)
\endCD
$$
where the map $\tau$
admits the following description:
An element of
$\times_{\ell}\Omega G$
is a
map $\phi$ from $\Delta_1$ to $G^{\ell}$ which sends the end points
to $e$;
now $\tau$
is given by the assignment to
$\phi$
of the composite
$$
\Delta_2
@>{(\phi \circ \eta^0,\phi \circ \eta^1)}>>
G^{\ell} \times G^{\ell}
@>r>>
G.
\tag4.3
$$
Inspection shows that
this composite indeed vanishes on the boundary
of $\Delta_2$ and hence passes to a based map
from the 2-sphere to $G$.
As already pointed out, the
left-hand vertical
fibration in (4.3) is rationally trivial,
cf. \cite\donalkro,
see also \cite\finite,
whence $\tau$ is rationally homotopically trivial.
However,
over the integers,
$\tau$ will in general {\it not\/}
be trivial.
For example,
let $Y$ be complex projective
2-space with the obvious cell decomposition,
so that $\ell = 1$ and the attaching map is the Hopf map
from $S^3$ to $S^2$.
Homotopically, the map $\tau$ then amounts to the map
from $\roman{Map}^o(S^2,BG)$
to $\roman{Map}^o(S^3,BG)$
induced by the Hopf map.
For simplicity, let $G=\roman {SU}(2)$.
Now $\roman \pi_4(\roman{Map}^o(S^2,BG))\cong \pi_6(BG)\cong \pi_5(S^3)$
which is cyclic of order 2, and we can represent the non-trivial element
by a non-trivial principal $\roman{SU}(2)$-bundle
on $S^4 \times S^2$.
In particular, its long exact homotopy sequence
will have non-trivial boundary operators.
Under the map from
$S^4 \times S^3$
to
$S^4 \times S^2$
induced by the Hopf map,
the bundle passes to
a principal $\roman{SU}(2)$-bundle
on $S^4 \times S^3$
having essentially the same long exact homotopy sequence
as the bundle
on $S^4 \times S^2$;
in particular, the bundle on
$S^4 \times S^3$
is non-trivial either.
A little thought reveals that this implies
that
the map
from $\roman{Map}^o(S^2,BG)$
to $\roman{Map}^o(S^3,BG)$
induced by the Hopf map is not null homotopic.

\medskip
\widestnumber\key{999}
\centerline{References}
\smallskip\noindent

\ref \no \anderone
\by D. W. Anderson
\paper A generalization of the Eilenberg-Moore spectral sequence
\jour Bull. Amer. Math. Soc.
\vol 78
\yr 1972
\pages  784--786
\endref

\ref \no \atibottw
\by M. Atiyah and R. Bott
\paper The Yang-Mills equations over Riemann surfaces
\jour Phil. Trans. R. Soc. London  A
\vol 308
\yr 1982
\pages  523--615
\endref

\ref \no \atiyjone
\by M. F. Atiyah and J. D. S. Jones
\paper Topological aspects of Yang-Mills theory
\jour Comm. Math. Phys.
\vol 61
\yr 1978
\pages  97--118
\endref

\ref \no \atiysemi
\by M. F. Atiyah, N. Hitchin, G. Segal, R. Lawrence
\paper Oxford seminar on Jones-Witten theory
\paperinfo Michaelmas Term
\yr 1988
\endref

\ref \no \bauesboo
\by H. J. Baues
\book Algebraic Homotopy
\publ Cambridge University Press
\publaddr Cambridge, England
\yr 1989
\endref

\ref \no \bottsega
\by R. Bott and G. Segal
\paper The cohomology of the vector fields on a manifold
\jour Topology
\vol 16
\yr 1977
\pages  285--298
\endref

\ref \no \bouskan
\by A. K. Bousfield and D. M. Kan
\paper Homotopy with respect to a ring
\jour Proc. Symp. Pure Math.
\vol 22
\yr 1971
\pages  59--64
\publ Amer. Math. Soc.
\publaddr Providence, R. I.
\endref

\ref \no \browhueb
\by R. Brown and J. Huebschmann
\paper Identities among relations,
\paperinfo
in: Low--dimensional topology, ed. R. Brown
and T. L. Thickstun
\jour
London Math. Soc. Lecture Note Series
\vol 48
\publ Cambridge Univ. Press, Cambridge, U. K.
\yr 1982
\pages 153--202
\endref

\ref \no \curtione
\by E. B. Curtis
\paper Simplicial homotopy theory
\jour Advances in Math.
\vol 6
\yr 1971
\pages 107--209
\endref

\ref \no \donalkro
\by S. K. Donaldson and P. B. Kronheimer
\book The geometry of four manifolds
\publ Oxford University Press
\publaddr Oxford, U. K.
\yr 1991
\endref

\ref \no \duponone
\by J. L. Dupont
\paper Simplicial de Rham cohomology and characteristic classes
of flat bundles
\jour Topology
\vol 15
\yr 1976
\pages  233--245
\endref

\ref \no \dwykanon
\by W. G. Dwyer and D. M. Kan
\paper Homotopy theory and simplicial groupoids
\jour Indag. Math.
\vol 46
\yr 1984
\pages 379--385
\endref

\ref \no \eilmactw
\by S. Eilenberg and S. Mac Lane
\paper On the groups ${\roman H(\pi,n)}$. I.
\jour Ann. of Math.
\vol 58
\yr 1953
\pages  55--106
\moreref
\paper II. Methods of computation
\jour Ann. of Math.
\vol 60
\yr 1954
\pages  49--139
\endref

\ref \no \afischon
\by A. E. Fischer
\paper A grand superspace for unified field theories
\jour General Relativity and Gravitation
\vol 18
\yr 1986
\pages 597--608
\endref

\ref \no \friehabe
\by T. Friedrich and L. Habermann
\paper The Yang-Mills equations on the two-dimensional sphere
\jour Comm. in Math. Phys.
\vol 100
\yr 1985
\pages 231--243
\endref

\ref \no \modus
\by J. Huebschmann
\paper Symplectic and Poisson structures of certain moduli spaces
\jour Duke Math. J. (to appear)
\paperinfo  hep-th 9312112
\endref

\ref \no \modustwo
\by J. Huebschmann
\paper Symplectic and Poisson structures of certain moduli spaces. II.
Projective representations of cocompact planar discrete groups
\jour Duke Math. J. (to appear)
\paperinfo  dg-ga/9412003
\endref

\ref \no \srni
\by J. Huebschmann
\paper Poisson geometry of certain moduli spaces
\paperinfo Lectures delivered at the 14th Winter school,
Czech Republic, Srni, January 1994
\jour Rendiconti del Circolo Matematico di Palermo, to appear.
\endref

\ref \no \finite
\by J. Huebschmann
\paper Extended moduli spaces and the Kan construction. II.
Lattice gauge theory
\paperinfo MPI preprint, 1995
\endref

\ref \no \huebjeff
\by J. Huebschmann and L. Jeffrey
\paper Group cohomology construction of symplectic forms
on certain moduli spaces
\jour Int. Math. Research Notices
\vol 6
\yr 1994
\pages 245--249
\endref

\ref \no \jamesboo
\by I. James
\book The topology of Stiefel manifolds
\bookinfo London Math. Soc. Lecture Note Series
\vol
\publ Cambridge University Press
\publaddr Cambridge, U. K.
\yr
\endref

\ref \no \jeffrtwo
\by L. Jeffrey
\paper Symplectic forms on moduli spaces
of flat connections on 2-manifolds
\paperinfo to appear in {\it Proceedings
of the Georgia International Topology Conference}, Athens, Ga.
1993, ed. by W. Kazez
\endref

\ref \no \jeffrthr
\by L. Jeffrey
\paper
Group cohomology construction of the cohomology of
moduli spaces
of flat connections
on 2-manifolds
\jour Duke Math. J. (to appear)
\endref

\ref \no \kanone
\by D. M. Kan
\paper On homotopy theory and c.s.s. groups
\jour Ann. of Math.
\vol 68
\yr 1958
\pages 38--53
\endref

\ref \no \kantwo
\by D. M. Kan
\paper A relation between CW-complexes and free c.s.s. groups
\jour Amer. J. of Math.
\yr 1959
\vol 81
\pages 512--528
\endref

\ref \no \kanthr
\by D. M. Kan
\paper A combinatorial definition of homotopy groups
\jour Ann. of Math.
\yr 1958
\vol 57
\pages 282--312
\endref

\ref \no \kanfou
\by D. M. Kan
\paper Homotopy groups, commutators, and $\Gamma$-groups
\jour Illinois J. of Math.
\yr 1960
\vol 4
\pages 1--8
\endref

\ref \no \karshone
\by Y. Karshon
\paper
An algebraic proof for the symplectic
structure of moduli space
\jour Proc. Amer. Math. Soc.
\vol 116
\yr 1992
\pages 591--605
\endref

\ref \no \kauflins
\by L. H. Kauffman and S. L. Lins
\book Temperley-Lieb Recoupling Theorey and Invariants of
3-manifolds
\bookinfo Annals of Mathematics studies,
No. 134
\publ Princeton University  Press
\publaddr Princeton, N. J.
\yr 1994
\endref

\ref \no \maclaboo
\by S. Mac Lane
\book Homology
\bookinfo Die Grundlehren der mathematischen Wissenschaften
 No. 114
\publ Springer
\publaddr Berlin $\cdot$ G\"ottingen $\cdot$ Heidelberg
\yr 1963
\endref

\ref \no \maclbotw
\by S. Mac Lane
\book Categories for the Working Mathematician
\bookinfo Graduate Texts in Mathematics
\vol 5
\publ Springer
\publaddr Berlin $\cdot$ G\"ottingen $\cdot$ Heidelberg
\yr 1971
\endref

\ref \no \maybook
\by P. J. May
\book Simplicial Objects in Algebraic Topology
\publ Van Nostrand
\yr 1967
\endref

\ref \no \milnothr
\by J. Milnor
\paper The geometric realization of a semi-simplicial complex
\jour Ann. of Math.
\vol 65
\yr 1957
\pages 357--362
\endref

\ref \no \philston
\by A. V. Phillips and D. A. Stone
\paper The computation of characteristic classes of lattice
gauge fields
\jour Comm. Math. Phys.
\vol 131
\yr 1990
\pages 255--282
\endref

\ref \no \pressega
\by A. Pressley and G. Segal
\book Loop Groups
\bookinfo Oxford Mathematical Monographs
\publ Oxford University Press
\publaddr Oxford
\yr 1986
\endref

\ref \no \puppeone
\by D. Puppe
\paper Homotopie und Homologie
in abelschen Gruppen und Monoidkomplexen. I. II
\jour Math. Z.
\vol 68
\yr 1958
\pages 367--406, 407--421
\endref

\ref \no \shulmone
\by H. B. Shulman
\book Characteristic classes and foliations
\bookinfo Ph. D. Thesis
\publ University of California
\yr 1972
\endref

\ref \no \steenone
\by N. Steenrod
\paper Milgram's classifying space of a topological group
\jour Topology
\vol 7
\yr 1968
\pages 349--368
\endref

\ref \no \weinstwe
\by A. Weinstein
\paper On the symplectic structure of moduli space
\paperinfo preprint 1992; A. Floer memorial, Birkh\"auser, to appear
\endref

\ref \no \whitehon
\by J. H. C. Whitehead
\paper On simply connected 4-dimensional polyhedra
\jour Comm. Math. Helv.
\vol 22
\yr 1949
\pages 48--92
\endref

\ref \no \whitehtw
\by J. H. C. Whitehead
\paper A certain exact sequence
\jour Ann. of Math.
\vol 52
\yr 1950
\pages 51--110
\endref

\ref \no \wittesix
\by E. Witten
\paper Quantum field theory and the Jones polynomial
\jour Comm. in Math. Phys.
\vol 121
\yr 1989
\pages 351--399
\endref

\enddocument